\documentclass{article}
\usepackage{amsfonts}

\usepackage{graphicx}
\usepackage{amsmath}

\newtheorem{theorem}{Theorem}

\newtheorem{lemma}[theorem]{Lemma}

\newenvironment{proof}[1][Proof]{\textbf{#1.} }{\ \rule{0.5em}{0.5em}}
\input{tcilatex}

\begin{document}

\begin{center}
{\huge \textbf{On a class of skew classical\\[10pt]
r-matrices with large carrier}}\\[20pt]

\textsc{\textbf{Vladimir D. Lyakhovsky} \\[0pt]
Theoretical Department, St. Petersburg State University,\\[0pt]
198904, St. Petersburg, Russia \\[30pt]
}
\end{center}

\begin{abstract}
The classical $r$-matrices with carriers $\frak{g}_c$ containing the Borel
subalgebra $\frak{b^+}(\frak{g})$ of simple Lie algebra $\frak{g}$ are
studied. Using the graphical presentation of the dual Lie algebra $\frak{g}%
^{\#}(r)$ we show that such solutions $r_{ech}$ of the CYBE always exist. To
obtain the explicit form of $r_{ech}$ we find the dual coordinates in which
the adjoint action of $\frak{g}_c$ can be reduced. This gives us the
detailed structure of the Jordanian $r$-matrices $r_{J}$ that are the
candidates for enlarging the initial full chain $r_{fch}$. We search the
desired solution $r_{ech}$ in the factorized form $r_{ech} \approx
r_{fch}+r_{J}$. This leads to the unique transformation: the canonical chain
is to be substituted with a special kind of peripheric $r$-matrices: $%
r_{fch}\longrightarrow r_{rfch}$. To illustrate the method the
case of $\frak{g}=sl(11)$ is considered in full details.
\end{abstract}

\section{Introduction}

Constant triangular solutions of the Yang-Baxter equation \cite{DR83}\ play
important role in algebra and applications
(for details see \cite{KR90,ESS00,AL} and references therein).
They describe Poisson structures $\pi
\left( r\right) $ compatible with the initial Lie algebra $\frak{g}$, i. e.
the mechanical systems that can exist on a space whose noncommutativity is
fixed by $\frak{g}$. Classical solutions $r$ define the nondegenerate bilinear
form $\omega _{r}$ on the space of a subalgebra $\frak{g}_{c}\subseteq \frak{%
g}$ called the carrier of $r$. Due to the properties of $\omega _{r}$ the
Poisson algebra $\frak{g}^{\#}\left( r\right) $ (induced by $r$) is
equivalent to $\frak{g}_{c}$. To describe all the structures $\pi \left(
r\right) $ for a given $g$ it is sufficient to find all the carrier
subalgebras $\frak{g}_{c}$ and the corresponding forms $\omega _{r}$. Such
subalgebras are called Lie Frobenius \cite{Ohn-Ela}.

For simple Lie algebras with the Cartan decomposition $\frak{g=n}_{-}+\frak{h%
}+\frak{n}_{+}=\frak{n}_{-}+\frak{b}_{+}=\frak{b}_{-}+\frak{n}_{+}$ there
are several classes of constant skew $r$-matrices that can be quantized in the
explicit form, i. e. for which the corresponding solutions $\mathcal{R}$ of
the quantum Yang-Baxter equation are known.

They are constructed in the form $\mathcal{R=F}_{21}\mathcal{F}^{-1}$ ,
where $\mathcal{F}$ is the twisting element, i.e. the solution of the twist
equations \cite{DR83}.

Among the twists for simple Lie algebras $\frak{g}$ the full
chains $\mathcal{F}_{ch}$ of extended Jordanian twists are distinguished by
the fact that their carriers $\frak{g}_{ch}$ might be large enough to
contain $\frak{n}_{+}$ (or $\frak{n}_{-}$). Chains $\mathcal{F}_{fch}$ (and
the corresponding solutions of CYBE $r_{fch}$ ) with $\frak{g}_{fch}\supset
\frak{n}_{\pm }$ are called full \cite{AKL}.

The properties of the full chains differ considerably for two main classes
of root systems $\Lambda ^{I}$ and $\Lambda ^{II}$. Consider the highest
root $\theta _{0}\in \Lambda \equiv \Lambda _{0}$ and the root subspace $%
V_{1}=\left\{ v\in V_{0}\mid \left( v,\theta _{0}\right) =0\right\} $ . In
the set $\Lambda _{1}=V_{1}\cap \Lambda _{0}$ find the highest root $\theta
_{1}$ and the subspace $V_{2}=\left\{ v\in V_{1}\mid \left( v,\theta
_{1}\right) =0\right\} $. This gives the sequence
\begin{equation*}
V_{0}\supset V_{1}\supset V_{2}\supset \ldots \supset V_{f}
\end{equation*}
that terminates when there are no roots in $V_{f}$ . (The last space $V_{f}$
is orthogonal to the roots $\theta _{0},\theta _{1},\ldots ,\theta _{f-1}$
.) For $\Lambda ^{I}$-systems the last space is nontrivial, $V_{f}\neq 0$,
while for $\Lambda ^{II}$ we have $V_{f}=0$.
When $\Lambda =\Lambda ^{II}$\ the
corresponding sequence of roots
\begin{equation*}
\theta _{0},\theta _{1},\ldots ,\theta _{f-1}
\end{equation*}
forms the basis of $V_{0}$\ and $f$ is equal to the rank of the algebra $%
\frak{g}$ . When $\Lambda =\Lambda ^{I}$\ we have $f<\mathrm{rank}\left(
\frak{g}\right) $ .

Series $B$, $C$ and $D_{t=2s}$ , are of the second type, $A_{l}$ series for $%
l>1$ and $D_{t=2s+1}$ are of the first. For $\frak{g}\left( \Lambda
^{II}\right) $ algebras the full chains of extended twists have the carrier $%
\frak{g}_{c}=\frak{b}_{+}\left( \frak{g}\right) $ . For $\frak{g}\left(
\Lambda ^{I}\right) $ algebras the carrier for the full chain contains $%
\frak{n}_{+}\left( \frak{g}\right) $ but do not contain the Cartan
subalgebra $\frak{h}\left( \frak{g}\right) $.

Most interesting is the situation for the series $A_{l}$ where the
dimension of $V_{f}$ increases with the rank.

We shall demonstrate that the full chains of $r$-matrices for algebras $%
\frak{g}\left( \Lambda ^{I}\right) $ can be enlarged by a series of
Jordanian summands. The resulting enlarged chain has the carrier $\frak{g}%
_{ech}$ that contains $\frak{b}_{+}\left( \frak{g}\right) $ and
nontrivially intersects with $\frak{n}_{-}$:
$\frak{g}_{ech}\cap \frak{n}_{-}\left( \frak{g}%
\right) =\frak{v}_{-}$. The space $\frak{v}_{-}$ in algebras $\frak{g}\left(
\Lambda ^{I}\right) $ is an analogue of the space $\frak{b}_{\perp }$ in
algebras $\frak{g}\left( \Lambda ^{II}\right) $ whose generators correspond
to the roots in $V_{f-1}$ (the space orthogonal to the set $\left\{ \theta
_{i}\mid i=1,...,f-1\right\} $). Consider an algebra $\frak{g}\left( \Lambda
^{II}\right) $ and the classical $r$-matrix for the full chain of extended
Jordanian terms $r_{\left( link\right) }$ : $r_{fch}=\sum_{i=1}^{f}r_{\left(
link\right) i}$ . Truncate the last term and consider the dual algebra $%
\frak{g}^{\#}\left( r_{\left( f-1\succ 1\right) }\right) $. For the
subalgebra $\frak{b}_{\perp }$ one can always find such a basis that all the
coproducts $\delta _{r_{\left( f-1\succ 1\right) }}\left( \frak{b}_{\perp
}\right) $ are primitive. We shall show that for the algebras of type I in
the full chain duals $\frak{g}^{\#}\left( r_{fch}\right) $
one can always find such a basis that all
the coproducts $\delta _{fch}\left( \frak{v}_{-}\right) $ are quasiprimitive. It
will be demonstrated that the quasiprimitivity in the case I (just as the
primitivity in the case II) permits to incorporate the additional Jordanian
terms in the corresponding chains $\left( r_{fch}\text{ and }r_{\left(
f-1\succ 1\right) }\right) $. In the second case this completes the full
chain $r_{fch}=r_{\left( f-1\succ 1\right) }+r_{J}$ with the carrier $\frak{g%
}_{fch}$ containing $\frak{h}\left( \frak{g}\right) $ . In the case $\frak{g}%
\left( \Lambda ^{I}\right) $ this leads to the enlarged chain $%
r_{ech}=r_{fch}+r_{J}$ with the similar property of the carrier: $\frak{g}%
_{ech}\supset \frak{h}\left( \frak{g}\right) $. The enlarged chains exist
for any algebra $\frak{g}\left( \Lambda ^{I}\right) $. We shall construct
them explicitly for $\frak{g}=sl\left( 2m+1\right) $.

\section{Properties of the generators $v_{-}$ in the dual algebra $\frak{%
g}^{\#}\left( r_{fch}\right) $}

From now on we shall consider the simple Lie algebras $\frak{g}$ of the
series $A_{n-1}$. The maximal sequence of highest roots $\left\{ \theta
_{s}\right\} _{s=0,\ldots ,f-1}$ has the length $f=\left\{
\begin{array}{c}
n/2 \\
\left( n-1\right) /2
\end{array}
\right. $ for even and odd $n$'s correspondingly. In particular we shall
concentrate our attention on the case $n=2m+1$ . The other cases can be
treated analogously. The canonical full chain of extended Jordanian terms
\cite{KLO} has the $r$-matrix
\begin{eqnarray}
r_{fch} &=&E_{\frac{n-1}{2},\frac{n+1}{2}}\wedge E_{\frac{n+1}{2},\frac{n+3}{%
2}}+H_{\frac{n-1}{2},\frac{n+3}{2}}\wedge E_{\frac{n-1}{2},\frac{n+3}{2}}+
 \notag\\
&&+\cdots +  \notag \\
&&+\sum_{k_{s+1}=s+2}^{n-s-1}E_{s+1,k_{s+1}}\wedge
E_{k_{s+1},n-s}+H_{s+1,n-s}\wedge E_{s+1,n-s}+  \notag \\
&&+\sum_{k_{s}=s+1}^{n-s}E_{s,k_{s}}\wedge E_{k_{s},n-s+1}\wedge
+H_{s,n-s+1}\wedge E_{s,n-s+1}+  \notag \\
&&+\sum_{k_{s-1}=s}^{n-s+1}E_{s-1,k_{s-1}}\wedge
E_{k_{s-1},n-s+2}+H_{s-1,n-s+2}\wedge E_{s-1,n-s+2}+  \notag \\
&&+\cdots + \notag\\
&&+\sum_{k_{1}=2}^{n-1}E_{1,k_{1}}\wedge E_{k_{1},n}+H_{1,n}\wedge E_{1,n}.
\label{ch}
\end{eqnarray}
This solution of the CYBE is a generic point of the $\left( m\right) $%
-dimensional variety of classical $r$-matrices
\begin{equation}
r_{fch}=\sum_{k=1}^{m}\xi _{k}\left( H_{k,n-k+1}\wedge
E_{k,n-k+1}+\sum_{s=k+1}^{n-k}E_{k,s}\wedge E_{s,n-k+1}\right) ,
\label{r-fch}
\end{equation}
where the independent parameters $\xi _{k}\in \mathsf{C}^{1}$ (one parameter
for each link of the chain; for details see (\cite{KLO})).

Consider the Lie bialgebra $\left( \frak{g},\frak{g}^{\#}\right) $
corresponding to the $r$-matrix (\ref{r-fch}). It was proved in \cite{L}
that for any simple Lie algebra $\frak{g}$ and the chain $r_{ch}$ the dual
algebra $\frak{g}^{\#}$ can be decomposed into the semidirect sum
\begin{equation}
\frak{g}^{\#}=\frak{g}_{c}\vdash \frak{a},  \label{dir-sum}
\end{equation}
with the Abelian ideal $\frak{a}$ generated by the elements $a_{j}^{\ast }$
dual to  $a_{j}\in \frak{g}\setminus \frak{g}_{c}.$ Let $\left\{
h_{\alpha },e_{\lambda }\mid \lambda \in \Lambda \left( \frak{g}\right)
,\alpha \in S\left( \frak{g}\right) \right\} $ be the Cartan-Weil basis of $%
\frak{g}$ with the root system $\Lambda \left( \frak{g}\right) $ and the set
of simple roots $S\left( \frak{g}\right) $. The elements $\left\{ h_{\alpha
}^{\ast },e_{\lambda }^{\ast }\right\} $ are canonically dual to $\left\{
h_{\alpha },e_{\lambda }\right\} $ and form a basis in $\frak{g}^{\#}$.
Algebra $\frak{g}^{\#}$ has the natural grading induced by the root system $%
\Lambda \left( \frak{g}\right) $ with the vector grading group $\Gamma
\left( r\right) $\ generated by certain subsets in $\Lambda \left( \frak{g}%
\right) $ (and the zero).

The selfdual subalgebra $\frak{g}_{c}\approx \left( \frak{g}_{c}\right)
^{\#} $ is graded according to the following rules \cite{L}:

\begin{itemize}
\item  each element $e_{\theta _{s}}^{\ast }$ has zero grade,

\item  each element $e_{h_{\theta _{s}}}^{\ast }$ has the grading vector $%
\left( -\theta _{s}\right) $,

\item  for each extension term $e_{\mu }\wedge e_{\nu } \in r_{ch}$ the
element $e_{\mu }^{\ast }$ has the grading $\left( -\nu \right) $ while $%
e_{\nu }^{\ast }$ has the grading $\left( -\mu \right) $.
\end{itemize}

Grading of the Abelian ideal $\frak{a}$ is inherited from the root system $%
\Lambda \left( \frak{g}\right) $ :

\begin{itemize}
\item  each element $e_{\xi }^{\ast }$ has the grading vector $\xi $,

\item  each element $h_{\perp }^{\ast }$ of the Cartan subalgebra $\frak{h}%
_{\perp }^{\#}$ has zero grade. (The corresponding element $h_{\perp }\in
\frak{h}$ is orthogonal to all the Cartan elements in the Jordanian summands
of the chain. When the chain is canonical the linear form induced by $%
h_{\perp }$ is orthogonal to all the roots $\left\{ \theta _{i}\mid
i=0,...,m-1\right\} $.)
\end{itemize}

These rules are valid for any type of chains including the peripheric ones
\cite{KwL}. Thus in the Jordanian summands $h_{\theta _{s}}\wedge e_{\theta
_{s}}$ the canonical Cartan factors $h_{\theta _{s}}$ can be changed for $%
h_{\theta _{s}}+\sum \gamma _{s}h_{s}^{\perp }$, $\gamma _{s}\in \mathsf{C}$%
. In the latter case the form induced on the grading space by the Cartan
element $e_{\theta _{s}}^{\ast }\in \left( \frak{g}_{c}\right) ^{\#}$ also
depends on $\gamma $'s. In particular the operator $ad_{e_{\eta }^{\ast }}$
with $e_{\eta }\in \frak{g}_{c}$ acts nontrivially on $h_{s}^{\perp \ast }$
iff $\eta \left( h_{s}^{\perp }\right) \neq 0$.

In the case of $r_{fch}$ the negative root vectors $\Lambda ^{-}\left( \frak{%
g}\right) $ (and the zero) constitute the (blue) set of grading vectors for
the generators of $\left( \frak{g}_{c}\right) ^{\#}$. The red set $\Lambda
^{-}\left( \frak{g}\right) \cup 0$ is attributed to the generators of $\frak{%
a}$ . Thus the compositions in $\frak{g}^{\#}$ are described by the
gradation diagram (we denote it by the same symbol $\Gamma \left( r\right) $%
) containing two superposed (coloured) sets: the algebra depicted by the
blue set is acting on the representation space drawn by the set of red
vectors. In our case these two sets are: $\left( \Lambda ^{-}\left( \frak{g}%
\right) \cup 0\right) ^{b}\cup \left( \Lambda ^{-}\left( \frak{g}\right)
\cup 0\right) ^{r}=\Gamma ^{b}\left( r\right) \cup \Gamma ^{r}\left(
r\right) =\Gamma \left( r\right) $.

Our task is to construct new
solutions for CYBE by adding to $r_{fch}$ the Jordanian type summands
\begin{equation*}
h_{\lambda }\wedge e_{\lambda },\qquad \lambda \in \Lambda ^{-}\left( \frak{g%
}\right) .
\end{equation*}
Consequently we are to fix all the
quasiprimitive generators in $\left( \frak{g},\frak{g}%
^{\#}\right) $. For this purpose we shall use the
gradation diagram $\Gamma \left( r\right) $. The
element $e_{\lambda }\in \frak{g}$ is quasiprimitive in the Lie bialgebra $%
\left( \frak{g},\frak{g}^{\#}\right)$ iff the grading
vector of its dual $e_{\lambda }^{\ast }$ cannot be presented as a sum of
two vectors one of which is nonzero blue. For a quasiprimitive $e_{\lambda
}$ the only adjoint operators whose image in the first derivative subalgebra
$\left( \frak{g}^{\#}\right) ^{\left( 1\right) }$ is $e_{\lambda }^{\ast }$ are
of the type $ad_{e_{\theta }^{\ast }}$ (with $e_{\theta }^{\ast }\in \frak{g}%
_{c}$). In other words we are to consider the end points in $\Gamma \left(
r\right) $ that cannot be reached from other red points by nonzero vectors
belonging to $\Gamma \left( r\right) ^{b}$.

Now we shall demonstrate that in $\left( \frak{g},\frak{g}^{\#}
 \right) $ we can have three sets of quasiprimitive
generators. The first set is $\left\{ e_{\theta _{s}}^{\ast }\mid s=1,\ldots
,m\right\} $.
(For the algebras of series $A_n$ we use the natural coordinate presentation
of the roots in which $\theta _{s}=e_{s}-e_{n-s+1}$.) These are the Cartan
elements of the subalgebra $\left( \frak{g}_{c}\right) ^{\#}$ , they have
(blue) zero grading vectors and are primitive. The second set contains the
elements $h_{\perp }^{\ast }\in \left( \frak{h}_{\perp }\right) ^{\#}$ dual
to the Cartan's in $\frak{h}_{\perp }$ . In the canonical chain the Cartan's
$h_{\perp }$
are orthogonal to all $\left\{ \theta _{k}\mid k=1,\ldots ,m\right\} $ ,
i.e. $\theta _{k}\left( h_{\perp }\right) =0$. These elements belong to $%
\frak{a}$, have (red) zero grading vectors and thus are primitive. The most
important for us is the third set -- $\left\{ e_{-\alpha }^{\#}\mid \alpha
\in S\left( \frak{g}\right) \right\} $. It contains the duals to negative
simple root generators $e_{-\alpha }$. These elements belong to $\frak{a}$
and have red grading vectors $\left\{ -\alpha _{k}=e_{k+1}-e_{k}\mid
k=1,\ldots ,2m\right\} $.
For the case of canonical chain (see the $r$-matrix (%
\ref{r-fch})) this set is empty. For any red vector
$\left( -\alpha _{k}\right) ^{\left( r\right) }$ we have the
corresponding blue one, $\left( -\alpha _{k}\right) ^{\left( b\right) }$,
that shifts the red zero point to the end point of $\left( -\alpha
_{k}\right) ^{\left( r\right) }$. Red zeros represent the elements of $%
\left( \frak{h}_{\perp }\right) ^{\#}$ and among them we can find those
whose duals are not orthogonal to $\alpha _{k}$. The adjoint operators
corresponding to $\left( -\alpha _{k}\right) ^{\left( b\right) }$ act
nontrivially on such elements $h_{\perp }^{\#}$ :
\begin{equation}
\left[ h_{\perp }^{\#},e_{\left( -\alpha _{k}\right) ^{\left( b\right)
}}^{\#}\right] ^{\#}\sim \left( -\alpha _{k}\left( h_{\perp }\right) \right)
e_{\left( -\alpha _{k}\right) ^{\left( r\right) }}^{\#}.
\label{n-quasi-terms}
\end{equation}
Such compositions are the only ones that violate the quasiprimitivity of $%
e_{\left( -\alpha _{k}\right) ^{\left( r\right) }}^{\#}$. We remind that the
subalgebra $\left( \frak{h}_{\perp }\right) ^{\#}$ depends on the choice of
Cartan generators in the Jordanian terms of the chain. For the canonical
chain they are $\left\{ h_{\theta _{s}}\right\} $. In the general case the
Jordanian terms look like $\left\{ \left( h_{\theta _{s}}+\eta
_{s}h_{s}^{\perp }\right) \wedge e_{\theta _{s}}\right\} $ . This means that
the sum $\sum \eta _{s}h_{s}^{\perp }\wedge e_{\theta _{s}}$ can be added to
the initial $r_{fch}$ :
\begin{equation*}
r_{rch}=r_{fch}+r_{r}=r_{fch}+\sum \eta _{s}h_{s}^{\perp }\wedge e_{\theta
_{s}}.
\end{equation*}
Notice that here $r_{r}$ is itself a solution of CYBE. The transformation $%
r_{fch}\Rightarrow r_{fch}+r_{r}$ signifies the rotation of the hyperplane $%
h_{\perp }^{\#}$ in the space $h^{\#}$. We call $r_{r}$ the \emph{rotation
term}.

Choosing the appropriate coefficients $\eta _{s}$ we can trivialize the
composition (\ref{n-quasi-terms}) and make the element $e_{\left( -\alpha
_{k}\right) ^{\left( r\right) }}^{\#}$ quasiprimitive. This means that the
space $h_{\perp }^{\#}$ becomes orthogonal to $\alpha _{k}$. In terms of $%
r_{rch}$ this means that the subspace $V_{ch}^{C}$ (generated by the duals
$\left\{ h_{s}^{\ast }\right\} $ of Cartan
factors belonging to the Jordanian terms $%
\left\{ h_{_{s}}\wedge e_{\theta _{s}}\mid s=1,\ldots ,m\right\} $) contains
$\alpha _{k}$. As far as \textrm{dim}$V_{ch}^{C}=m$ such requirements can be
fulfilled at most for $m$ simple roots $\alpha _{k}\in S_{ch}\left( \frak{g}%
\right) $. Below we shall show that for $\frak{g}=sl(n)$ these $m$
additional conditions can be always fulfilled and we can get exactly $m$
quasiprimitive generators $\left\{ e_{-\alpha _{k}}^{\#}\mid \alpha _{k}\in
S_{ch}\left( \frak{g}\right) \right\} $.

\section{Construction of the enlarged chain}

To construct the additional jordanian terms $%
\left\{ h_{\lambda }\wedge e_{\lambda },\mid \lambda \in \Lambda ^{-}\left(
\frak{g}\right) \right\} $ we need the subset of
mutually orthogonal negative simple roots. The reason is that any nontrivial
combination of the Jordanian terms $\sum h_{_{k}}^{\perp }\wedge e_{-\alpha
_{k}}$ with noncommuting generators $e_{-\alpha _{k}}$ violates CYBE and the
nonzero Schouten bracket $\left[ \left[ \sum h_{_{k}}^{\perp }\wedge
e_{-\alpha _{k}}\right] \right] $ cannot be compensated by the terms like $%
\left[ \left[ r_{rch},\sum h_{_{k}}^{\perp }\wedge e_{-\alpha _{k}}\right] %
\right] +\left[ \left[ \sum h_{_{k}}^{\perp }\wedge e_{-\alpha _{k}},r_{rch}%
\right] \right] $. We shall consider the following set of $m$
mutually orthogonal roots
\begin{equation}
S_{ch}^{\perp }\left( \frak{g}\right) =\left\{ \alpha
_{k}=e_{2k-1}-e_{2k}\mid k=1,\ldots ,m\right\} .  \label{ort-roots}
\end{equation}

The adjoint action (\ref{n-quasi-terms}) is trivialized iff
\begin{equation}
\mathcal{L}\left( h_{s}^{\ast } \mid s=1,\ldots ,m\right)
\supseteq S_{ch}^{\perp
}\left( \frak{g}\right) ,  \label{ort-cond}
\end{equation}
(Here $\mathcal{L}\left( h_{s}^{\ast }\right) $ is the space generated by $%
\left\{ h_{s}^{\ast }\right\} $ .)

\begin{lemma}
The $r$-matrix
\begin{eqnarray*}
r_{r} &=&\sum_{i=1}^{m}\xi _{i}\widetilde{H}_{i}^{\perp }\wedge E_{i,n-i+1},
\\
\widetilde{H}_{i}^{\perp } &=&\sum_{j=i}^{n-i}\left( -1\right)
^{j+1}H_{j,j+1}
\end{eqnarray*}
rotates the chain
\begin{equation*}
r_{fch}=\sum_{k=1}^{m}\xi _{k}\left( H_{k,n-k+1}\wedge
E_{k,n-k+1}+\sum_{p=k+1}^{n-k}E_{k,p}\wedge E_{p,n-k+1}\right) .
\end{equation*}
The sum
\begin{equation*}
r_{rch}=r_{fch}+r_{r}
\end{equation*}
is the solution to CYBE and in the dual algebra $\frak{g}%
^{\#}\left( r_{rch}\right) $ the elements \\
$\left\{ e_{-\alpha _{k}}^{\#}\mid
\alpha _{2k-1}=e_{2k-1}-e_{2k}\mid 1,\ldots ,m\right\} $ are quasiprimitive.
\end{lemma}

\begin{proof}
Elements $\widetilde{H}_{i}^{\perp }$ are orthogonal to the highest roots $%
\left\{ \theta _{s}\right\} _{s=1,\ldots ,m}$ of the chain (\ref{ch}). In
the algebra $\frak{g}^{\#}\left( r_{fch}\right) $ we have the set of
primitive commuting generators: $\left\{ \left( E_{i,n-i+1}\right) ^{\ast
},\left( \widetilde{H}_{i}^{\perp }\right) ^{\ast }\mid i,j=1,\ldots
,m\right\} $. It follows that both $r_{r}$ and $r_{fch}+r_{r}$\ are the
solutions to the CYBE. The quasiprimitivity is realized for the basic
elements of $\frak{g}^{\#}$ correlated with the decomposition (\ref{dir-sum}%
). For the set $\left\{ e_{-\alpha }^{\#}\mid \alpha \in S\left( \frak{g}%
\right) \right\} $ this means that in terms of the duals $e_{-\alpha
_{k}}^{\ast }$the generators $e_{-\alpha }^{\#}$ have the following
expressions:
\begin{eqnarray*}
e_{-\alpha _{t}}^{\#} &=&\left( e_{-\alpha _{t}}^{\ast }+e_{\alpha
_{n-t}}^{\ast }\right) \qquad t=1,\ldots ,q-1,\widehat{q},q+1\ldots ,n-1; \\
e_{-\alpha _{q}}^{\#} &=&e_{-\alpha _{q}}^{\ast }.\qquad q=\left\{
\begin{array}{c}
m\text{ for odd }m \\
m+1\text{ for even }m
\end{array}
\right\}
\end{eqnarray*}
Rewriting the $r$-matrix
\begin{eqnarray*}
r_{rch} &=&r_{fch}+r_{r}= \\
&=&\sum_{k=1}^{m}\xi _{k}\left( \left( H_{k,n-k+1}+H_{k}^{\perp }\right)
\wedge E_{k,n-k+1}+\sum_{p=k+1}^{n-k}E_{k,p}\wedge E_{p,n-k+1}\right)
\end{eqnarray*}
and the Cartan elements
\begin{eqnarray*}
&&\left\{ H_{k,n-k+1}+\widetilde{H}_{k}^{\perp }\mid k=1,\ldots ,m\right\}
\\
&=&\left\{
\begin{array}{c}
\left( H_{1,2}+\ldots +H_{2m-1,2m}\right) ,\left( H_{3,4}+\ldots
+H_{2m-1,2m}\right) , \\
\left( H_{3,4}+\ldots +H_{2m-3,2m-2}\right) ,\ldots ,H_{m,m+1}
\end{array}
\right\} .
\end{eqnarray*}
we notice that the linear span $\mathcal{L}\left\{ H_{k,n-k+1}+\widetilde{H}%
_{k}^{\perp }\mid k=1,\ldots ,m\right\} $ contains the set $\left\{
H_{2k-1,2k}\quad \mid \quad k=1,\ldots ,m\right\} $ of Cartan elements dual
to the simple roots $\left\{ \alpha _{2k-1}=e_{2k-1}-e_{2k}\mid 1,\ldots
,m\right\} $. The subspace $\left( \frak{h}_{\perp }\right) ^{\#}$ in $%
\frak{g}^{\#}\left( r_{rch}\right) $ is defined as orthogonal to the set of
Jordanian Cartan generators in $r_{rch}$. As we have seen this is
equivalent to the requirement that $\frak{h}_{\perp }^{\ast }$ is orthogonal
to the roots $\left\{ \alpha _{2k-1}\mid 1,\ldots ,m\right\} $ . As a result
the operators $\mathrm{ad}_{\left( e_{\left( -\alpha _{k}\right) ^{\left(
b\right) }}^{\#}\right) }$ cannot shift the elements of $\left( \frak{h}%
_{\perp }\right) ^{\#}$ in the direction
$\left(-\alpha _{k}\right)^{\left( r\right) }$ and the elements
$e_{
\left(-\alpha _{k}\right)^{\left( r\right) }}^{\#}$ remain quasiprimitive.
\end{proof}

We have proved that there is a possibility to enlarge $r_{rch}$ by the set
of $m$ Jordanian terms formed by independent Cartan elements $\left\{
\widetilde{H}_{i}^{\perp }\mid j=1,\ldots ,m\right\} $ (belonging to the
subspace of primitive generators $\left( \frak{h}_{\perp }\right) ^{\#}$
in the dual algebra $\frak{g}^{\#}\left( r_{rch}\right) $) ) and the
elements
\begin{equation}
\left\{
\begin{array}{c}
\left( e_{-\alpha _{2k-1}}+e_{\alpha _{n-2k+1}}\right) ,e_{-\alpha
_{2p-1}}\mid k=1,\ldots ,p-1,\widehat{p},p+1,\ldots ,m \\
p=\left\{
\begin{array}{c}
\frac{m+1}{2}\text{ for odd }m \\
\frac{m+2}{2}\text{ for even }m
\end{array}
\right\}
\end{array}
\right\} ,  \label{jord-coord}
\end{equation}
corresponding to the set of quasiprimitive generators $\left\{ e_{\left(
-\alpha _{2k-1}^{\left( r\right) }\right) }^{\#}\!\mid k=1,\ldots ,m\right\}
$. The explicit form for such enlargement is to be found now.

In the case of $sl(n)$ , $n=2m+1$ the set (\ref{jord-coord}) can be written
in terms of matrix units $E_{i,j}$ :
\begin{equation}
\left\{
\begin{array}{c}
\widehat{E}_{k}=E_{2k,2k-1}+E_{n-2k+1,n-2k+2},\mid k=1,\ldots ,p-1,\widehat{p%
},p+1,\ldots ,m \\
\widehat{E}_{p}=E_{2p,2p-1};
\end{array}
\right\} .  \label{jord-coord-mat}
\end{equation}
To construct $m$ independent Jordanian terms we must find for each $\widehat{%
E}_{k}$ the corresponding $H_{k}^{\perp }$ with the properties:
\begin{eqnarray}
\left[ H_{k}^{\perp },\widehat{E}_{l}\right] &=&\delta _{kl}\widehat{E}%
_{k},\qquad k,l=1,\ldots ,m;  \label{orth-cart-cond} \\
H_{k}^{\perp } &\in &\frak{h}_{\perp }.  \notag
\end{eqnarray}
This gives $2m+1$ relations that completely fix the set $\left\{
H_{k}^{\perp }\right\} $ :
\begin{eqnarray}
H_{k}^{\perp } &=&\sum_{v=1}^{n}\frac{4k-2}{n}E_{v,v}-\sum_{u=1}^{2k-1}%
\left( E_{u,u}+E_{n-u+1,n-u+1}\right) ,  \label{cart-org-jord} \\
k &=&1,\ldots ,m.  \notag
\end{eqnarray}
Thus we have found the exact form of the Jordanian terms $\sum h_{\perp
_{k}}\wedge e_{-\alpha _{k}}$ :
\begin{equation}
r_{J}=\sum_{k=1}^{m}H_{k}^{\perp }\wedge \widehat{E}_{k}.  \label{jord-terms}
\end{equation}
Notice that this expression itself is obviously the solution to CYBE, that
is the pair $\left( \frak{g},\frak{g}^{\#}\left( r_{J}\right) \right) $ is a
Lie bialgebra as well as $\left( \frak{g},\frak{g}^{\#}\left( r_{rch}\right)
\right) $. The question is whether the pair $\left( \frak{g},\frak{g}%
^{\#}\left( r_{rch}+r_{J}\right) \right) $ is a Lie bialgebra or in other
words whether the sum of the $r$-matrices is again the solution to CYBE.

One simple example shows that the latter is not true. In the almost trivial
and well studied case $m=1$ ($n=3$) it was demonstrated (\cite{L},
\cite{LSam}) that the sum $r_{rch\left( 3\right) }+r_{J\left( 3\right) }$ does
not satisfy the CYBE. We'll see below that the same is true for any algebra $%
sl(n)$. It was found out (\cite{L}, \cite{LSam})  that the $r_{rch\left(
3\right) }$-matrix
\begin{equation*}
r_{rch\left( 3\right) }=2H_{12}\wedge E_{13}+E_{12}\wedge E_{23}
\end{equation*}
can be enlarged by a Jordanian like summand only if the latter is deformed:
\begin{equation}
r_{dj}=H^{\perp }\wedge \left( E_{21}+2E_{13}\right) ,\qquad H^{\perp }=-%
\frac{1}{3}\left( E_{11}-2E_{22}+E_{33}\right) .  \label{def-jord}
\end{equation}
Notice that according to the structure (\ref{jord-coord-mat}) of the dual
coordinates $\widehat{E}%
_{k}$  in the case $\frak{g}=sl(3)$ we have $\widehat{E%
}_{1}=E_{21}$ and the deformed expression (\ref{def-jord}) differs from the
additional Jordanian term $r_{J\left( 3\right) }=H_{1}^{\perp }\wedge
\widehat{E}_{1}$.The deformation in $r_{dj}$ results in and is coordinated
with the quasiprimitivity of $e_{\left( -\alpha _{1}\right) ^{\left( r\right)
}}^{\#}$ in $\left( sl(3),sl(3)^{\#}\left( r_{rch}\right) \right) $:
\begin{equation*}
\delta _{rch}\left( E_{21}\right) =2E_{13}\wedge E_{21}.
\end{equation*}
\qquad This example shows that we can redefine the Jordanian term in the
chain
\begin{equation*}
r_{rch\left( 3\right) }^{\prime }=2\left( H_{12}+H^{\perp }\right) \wedge
E_{13}+E_{12}\wedge E_{23},
\end{equation*}
and then consider the $r$-matrix
\begin{equation*}
r_{J}=H^{\perp }\wedge E_{21}
\end{equation*}
as an additional Jordanian for the twice rotated chain $r_{rch\left(
3\right) }^{\prime }$,
\begin{equation*}
r_{rch\left( 3\right) }^{\prime }+r_{J}=r_{rch\left( 3\right) }+r_{dj}.
\end{equation*}
To obtain the necessary result here two conditions are to be satisfied: $%
\widehat{E}$ is to be an element of the inverse root in one of the extension
terms ($E_{-\alpha _{1}}$ in our case) and the Cartan element in the
corresponding link of the chain (here $\widehat{H}=2\left( H_{12}+H^{\perp
}\right) $ ) must obey the relation $\widehat{H}-\left[ E_{\alpha
_{1}},E_{-\alpha _{1}}\right] \sim H^{\perp }$.

We start considering the general case supposing that the similar scheme
works for $\frak{g}=sl(n)$ and in particular for $n=2m+1$ .

First let us remind that in an arbitrary the full chain
\begin{equation}
r_{rch}=\sum_{k=1}^{m}\left( \widehat{H}_{k}\wedge
E_{k,n-k+1}+\sum_{p=k+1}^{n-k}E_{k,p}\wedge E_{p,n-k+1}\right) ,
\label{ch-modified}
\end{equation}
the Cartan elements are subject to the conditions
\begin{equation*}
\left[ \widehat{H}_{k},E_{l,n-l+1}\right] =\delta _{kl}E_{l,n-l+1}.
\end{equation*}
This is equivalent to the requirement
\begin{equation}
\widehat{H}_{k}=H_{k,n-k+1}+\sum \gamma _{l}H_{l}^{\perp }.
\label{extrem-cart}
\end{equation}
Take the additional Jordanian terms as in (\ref{jord-terms}) and consider
the sum
\begin{equation}
r_{ech}=r_{rch}+r_{J}=\sum_{l=1}^{m}\left( \widehat{H}_{l}\wedge
E_{l,n-l+1}+\sum_{p=l+1}^{n-l}E_{l,p}\wedge E_{p,n-l+1}+H_{l}^{\perp }\wedge
\widehat{E}_{l}\right) .  \label{extrem-r}
\end{equation}
To find the expressions for the modified Cartan's $\widehat{H}_{l,n-l+1}$
impose the CYBE,
\begin{equation}
\left[ \left[ r_{ech}\right] \right] =\left[ \left[ r_{rch},r_{J}\right] %
\right] +\left[ \left[ r_{J},r_{rch}\right] \right] =0.  \label{extrem-cybe}
\end{equation}
This leads to the following relations (in the first set we supplement the
definition of $\widehat{H}_{k}$ \ by putting $\widehat{H}_{m+1}=0$):

\begin{equation*}
\widehat{H}_{k}-\widehat{H}_{k+1}-2H_{\left( 2\chi \left( k\right) +1\right)
,2\chi \left( k\right) }\sim H_{\chi \left( k\right) }^{\perp },\quad \left.
\begin{array}{c}
k=1,\ldots ,m, \\
\chi \left( k\right) =\frac{1}{4}\left( n+2+\left( n-2k\right) \left(
-1\right) ^{k}\right) ;
\end{array}
\right.
\end{equation*}
\begin{equation*}
\sum_{k}^{m}\left[ \widehat{H}_{k},E_{2l+2,2l+1}\right] =0,\qquad l=1,\ldots
,m;
\end{equation*}
\begin{equation*}
\left.
\begin{array}{l}
\left[ \widehat{H}_{s},E_{\psi \left( s\right) ,\psi \left( s\right) +1}%
\right] =E_{\psi \left( s\right) ,\psi \left( s\right) +1}, \\
\left[ \widehat{H}_{s+1},E_{\psi \left( s\right) ,\psi \left( s\right) +1}%
\right] =-E_{\psi \left( s\right) ,\psi \left( s\right) +1},
\end{array}
\right\} ~\left.
\begin{array}{c}
s=1,\ldots ,m-1, \\
\psi \left( s\right) =\frac{1}{2}\left( n+\left( n-2s\right) \left(
-1\right) ^{s+1}\right) ;
\end{array}
\right.
\end{equation*}
Among these relations independent are only the first two sets. They have the
unique solution. Using the explicit expression (\ref{cart-org-jord}) for $%
H_{q}^{\perp }$  we get the final answer:
\begin{equation}
\widehat{H}_{k}=\left( -1\right) ^{k+1}\left(
\begin{array}{c}
\frac{2k-1}{n}\sum_{i=1}^{n}E_{i,i}-\sum_{j=1}^{k-1}\left(
E_{j,j}+E_{n-j+1,n-j+1}\right) \\
+\frac{1}{2}\left(
\begin{array}{c}
\left( -1\right) ^{k+1}\left( E_{k,k}-E_{n-k+1,n-k+1}\right) \\
-E_{k,k}-E_{n-k+1,n-k+1}
\end{array}
\right)
\end{array}
\right) .  \label{cart-chain}
\end{equation}
When the Cartan elements (\ref{extrem-cart}) are chosen in the form
(\ref{cart-chain}) the $r$-matrix $r_{ech}$ is the solution of CYBE with
the carrier $\frak{g}_{ech}$ containing $\frak{b}_{+}$ and $m$ elements
$\widehat{E}_{k}$ from the negative sector of Cartan decomposition. These
elements correspond to the subset of $m$ mutually
orthogonal negative simple roots $\left\{ -\alpha _{k}\mid \alpha _{k}\in
S_{ch}^{\perp }\left( \frak{g}\right) \right\} $ (\ref{ort-roots}).
Thus we have proved the following statement:

\begin{theorem}
The full chain of extended $r$-matrices for $sl(2m+1)$ can be uniquely
enlarged by the set of $m$ independent Jordanian terms.
\end{theorem}

\section{Parameterization}

It is well known that in fact the canonical $r$-matrices are the
representatives of the parameterized sets of objects. Such parameterization
is induced by the automorphisms of the carrier subalgebra of the $r$-matrix
(see \cite{KLO} for the case of full chain). In our case the carrier $\frak{g%
}_{ech}$ contain $\frak{b}_{+}$ and $m$ elements $\widehat{E}_{k}$ belonging
to $\frak{n}_{-}$.
The following automorphism of $\frak{g}_{ech}$ leads to the parameterized
variety of $r$-matrices:

\begin{itemize}
\item  for the chain carrier $\frak{g}_{fch}$\ :
\end{itemize}

For each link with the number $i$ we perform the transformations
\begin{eqnarray*}
E_{a,b} &\Longrightarrow &\xi _{i}E_{a,b} \\
E_{c,d} &\Longrightarrow &\frac{1}{\xi _{i}}E_{c,d}
\end{eqnarray*}
where the indices for odd $i$'s are
\begin{eqnarray}
a &=&i,\ldots ,n-i;\quad b=n+1-i,\ldots ,n;  \label{var-1} \\
c &=&1,\ldots ,i-1;\quad d=i,\ldots ,n-i;  \notag
\end{eqnarray}
and for even $i$'s --
\begin{eqnarray}
a &=&1,\ldots ,i;\quad b=i+1,\ldots ,n+1-i;  \label{var-2} \\
c &=&i+1,\ldots ,n+1-i;\quad d=n+2-i,\ldots ,n.  \notag
\end{eqnarray}
(For the first link the zone for $\frac{1}{\xi _{1}}$ degenerates to zero.)

\begin{itemize}
\item  for the quasijordanian part the parameterization is induced by the
transformations:
\end{itemize}

\begin{eqnarray*}
E_{a,b} &\Longrightarrow &\zeta _{i}E_{a,b}; \\
E_{c,d} &\Longrightarrow &\frac{1}{\zeta _{i}}E_{c,d};
\end{eqnarray*}
where for $i\leq \frac{m+1}{2}$ we have
\begin{eqnarray}
a &=&2i,\quad b=2i-1; \\
a &=&2i,\ldots ,n-2i+1;\quad b=n-2i+2,\ldots ,n;  \notag \\
c &=&1,\ldots ,2i-1;\quad d=2i,\ldots ,n-2i+1.  \notag
\end{eqnarray}
and for $i>\frac{m+1}{2}$ --
\begin{eqnarray}
a &=&2i,\quad b=2i-1; \\
a &=&1,\ldots ,n-2i+1;\quad b=n-2i+2,\ldots ,2i-1;  \notag \\
c &=&n-2i+2,\ldots ,2i-1;\quad d=2i,\ldots ,n.  \notag
\end{eqnarray}

Applying the above automorphism we get
\begin{equation}
\begin{array}{l}
r_{ech}\left( \xi ,\zeta \right) =r_{rch}\left( \xi \right) +r_{J}\left( \xi
,\zeta \right) = \label{r-ext-param} \\[3mm]
= \sum_{l=1}^{m}\left( \xi _{l}\widehat{H}_{l}\wedge E_{l,n-l+1}+\xi
_{l}\sum_{p=l+1}^{n-l}E_{l,p}\wedge E_{p,n-l+1}+\zeta _{l}H_{l}^{\perp
}\wedge \widehat{E}_{l}\left( \xi \right) \right) . \notag
\end{array}
\end{equation}
Here the generators $\widehat{E}_{k}\left( \xi \right) $ depend on the
parameters of the chain:
\begin{eqnarray}
\widehat{E}_{k}\left( \xi \right) &=&E_{2k,2k-1}+\frac{\xi _{2k-1}}{\xi _{2k}%
}E_{n-2k+1,n-2k+2},\qquad k=1,\ldots ,p-1;  \label{coord-param} \\
\widehat{E}_{p}\left( \xi \right) &=&E_{2p,2p-1};  \notag \\
\widehat{E}_{k}\left( \xi \right) &=&E_{2k,2k-1}+\frac{\xi _{n-2k+1}}{\xi
_{n-2k+2}}E_{n-2k+1,n-2k+2},\qquad k=p+1,\ldots ,m.  \notag
\end{eqnarray}
Notice that in the variety $r_{ech}\left( \xi ,\zeta \right) $ the
parameters $\zeta _{k}$ and $\xi _{1}$ are in $C^{1}$ while $\left\{ \xi
_{l}\mid l=2,\ldots ,m\right\} \in C^{1}\setminus 0$ . We cannot arbitrarily
switch off the links of the chain in $r_{exch}\left( \xi ,\zeta \right) $.
This can be done successively from the first to the $m$-th or in all the chain
simultaneously. Putting all $\xi _{k}$ 's equal zero we get the set of
purely Jordanian terms where the second summands in $\widehat{E}_{k}\left(
\xi \right) $ have indeterminate coefficients. Obviously the additional
Jordanians can be switched off in any order. \ When all the $\zeta _{k}$ 's
are zeros the full chain is restored with the special choice (\ref
{cart-chain}) of the Cartan's $\widehat{H}_{l}$.

\section{Conclusions}

We have proved that the full chain $r$-matrices for $sl(2m+1)$ can be
enlarged by the sum of independent Jordanian terms. The number of such terms
can be equal to $\mathrm{rank}\left( \frak{g}\right) -m$ so that the carrier
$\frak{g}_{ech}$ of the enlarged chain $r_{ech}$ contains the Borel
subalgebra $\frak{b}_{+}\left( \frak{g}\right) $ and $m$ commuting negative
simple root generators $\left\{ e_{\left( -\alpha _{2k-1}\right) }\mid
k=1,\ldots ,m\right\} $.

Similar result can be obtained for $sl(2m)$ algebras. Here the general
structure of the enlarged chain is the same
\begin{equation*}
r_{ech}=r_{rch}+r_{J}=\sum_{l=1}^{m}\left( \widehat{H}_{l}\wedge
E_{l,n-l+1}+\sum_{p=l+1}^{n-l}E_{l,p}\wedge E_{p,n-l+1}\right)
+\sum_{k=1}^{m-1}H_{k}^{\perp }\wedge \widehat{E}_{k}.
\end{equation*}
The difference (with the odd case considered above) is that here the
dimension of the space $\frak{h}_{\perp }$ is $2m-1-f=m-1$. This provides
additional freedom in choosing the set $S_{ch}^{\perp }\left( \frak{g}%
\right) $ of orthogonal simple roots. As far as for even $n$ the last link
of the full chain is degenerate (contains only the Jordanian term) there are
no such coordinates as $\widehat{E}_{p}$ with $e_{-\alpha
_{q}}^{\#}=e_{-\alpha _{q}}^{\ast }$. All this leads to different sets of $%
\left\{ H_{k}^{\perp }\right\} $ and $\left\{ \widehat{H}_{l}\right\} $.
Nevertheless the main property is valid: the enlarged chain $r_{ech}$ exists
in the form $r_{rch}+r_{J}$, its carrier contains $\frak{b}_{+}\left( \frak{g%
}\right) $ and $m-1$ commuting negative simple root generators.

As it is explained in the Appendix A among the classical Lie algebras there
exists only one more set of the type I -- the even-odd subset of orthogonal
algebras $D_{t=2s+1}$ . Using the same tools as in Section 3 one can
check that in this case the enlarged chains also exists. The space $\frak{h}%
_{\perp }$ is always one-dimensional and as a result the additional
Jordanian part $r_{J}$ contains only one term.

We come to the conclusion that for any simple Lie algebra there always
exists a chain of extended Jordanian $r$-matrices whose carrier contains $%
\frak{b}_{+}\left( \frak{g}\right) $. For algebras $\frak{g}\left( \Lambda
^{II}\right) $ these are the full chains $r_{fch}$ , for algebras $\frak{g}%
\left( \Lambda ^{II}\right) $ -- the enlarged chains $r_{ech}$.

For algebras of series $A_{n-1}$ the set of parabolic $r$-matrices $b_CG{}(n)$
was constructed in \cite{GG97}. The sets $b_{CG}(n)$ and $r_{ech}(n)$ intersects
at the point $n=3$: the enlarged chain carrier $\frak{g}_{ech}$ is the
parabolic subalgebra when $\frak{g}=\frak{sl}(3)$.

The next problem is to quantize explicitly the $r$-matrices $%
r_{ech}$ or to construct the solutions $\mathcal{F}_{ech}$ of the twist
equations \cite{DR83} corresponding to $r_{ech}$. It is obvious that such
twist as $\mathcal{F}_{ech}$ will contain the full chain twist $\mathcal{F}%
_{fch}$ as a factor, $\mathcal{F}_{ech}=\mathcal{F}_{J}\mathcal{F}_{rch}$ ,
thus it is sufficient to find the factor $\mathcal{F}_{J}$. It must be taken
into account that the form of $\mathcal{F}_{J}$ might differ from the canonical
exponent $e^{H\otimes E}$. It still can be deformed by the chain factor $%
\mathcal{F}_{rch}$. The results obtained above show that in the first order
(with respect to the overall deformation parameter) such deformation can be
eliminated, the appropriate coordinate transformation always exists. A
separate study will be devoted to the construction of $\mathcal{F}_{ech}$
and related topics.

\section{Acknowledgments}
This work was supported by Russian Foundation for Basic Research under the grant
N 03-01-00593.

\section{Appendix A. The sequence $V_{0}\supset V_{1}\supset
\ldots \supset V_{f-1}$ for simple Lie algebras}

\begin{itemize}
\item  Series $A$: for $sl(n)$ ($A_{n-1}$) the roots are
\begin{equation*}
e_{i}-e_{j}\qquad i,j=1,2,\ldots ,n.
\end{equation*}
\ Let $e_{1}-e_{n}$ be the highest root. The subsystems are
\begin{equation*}
\Lambda _{k}=\left\{ e_{i}-e_{j}\mid i,j=\left( k+1\right) ,\ldots ,\left(
n-k\right) ;\right\}
\end{equation*}
There is a difference between even and odd cases, $n=2m$ and $n=2m+1$. The
maximal length of the sequence of $\left\{ \theta _{s}\right\} $ is $f=m$ .
The dimension dim$V_{f-1}$ is
\begin{equation*}
\mathrm{dim}V_{f-1}=r-m=\left\{
\begin{array}{lll}
m-1 & \mathrm{for} & n=2m \\
m & \mathrm{for} & n=2m+1
\end{array}
\right.
\end{equation*}

\item  Series $B$: for $so(2m+1)$ ($B_{p}$) the roots are
\begin{equation*}
\left\{ \pm e_{i},\pm e_{i}\pm e_{j}\mid i,j=1,2,\ldots ,m;\right\}
\end{equation*}
The highest root is $e_{1}-e_{p}$. In the subspace $V_{1}$ we shall find the
roots
\begin{equation*}
\left\{
\begin{array}{c}
\pm e_{i},\pm e_{i}\pm e_{j}\mid i,j=2,\ldots ,m-1; \\
e_{1}+e_{m},\qquad -e_{1}-e_{m}.
\end{array}
\right\}
\end{equation*}
and $\mathrm{dim}V_{1}=m-1$. The maximal length is $m$ and is equal to the
rank. Algebras $B_{m}$ are of the type II.

\item  Series $C$: for $sp(m)$ ($C_{m}$) the roots are
\begin{equation*}
\left\{ \pm 2e_{i},\pm e_{i}\pm e_{j}\mid i,j=1,2,\ldots ,m;\right\}
\end{equation*}
The highest root is $2e_{1}$ and the sequence
\begin{equation*}
\theta _{0}=2e_{1},\theta _{1}=2e_{2},\ldots ,\theta _{p-1}=2e_{m}
\end{equation*}
contains $m$ elements and forms the basis of the root space $V_{0}$.
Algebras $C_{m}$ are of the type II.

\item  Series $D$: for $so(2m)$ ($D_{m}$) the roots are
\begin{equation*}
\left\{ \pm e_{i}\pm e_{j}\mid i\neq j;i,j=1,2,\ldots ,m;\right\}
\end{equation*}
Here the properties depend on $m$ . For $m$ even $D_{m}$ is of type II. For $%
m$ odd it is of type I. The latter is obvious due to the isomorphism $%
so(6)\approx sl(4)$. The situation with even-odd $so$-algebra differs from
the $A$-series case. Here the dimension\ of $V_{m}$ is always one: $\mathrm{%
dim}V_{m}=1.$
\end{itemize}

\section{Appendix B.}

Here we present the construction of the enlarged full chain for the typical
case of $\frak{g}=sl(11)$.

The full rotated chain $r$-matrix for $sl(11)$ is
\begin{equation*}
r_{rch}\left( \gamma \right) =\sum_{k=1}^{5}\left( \widehat{H}_{k}\left(
\gamma \right) \wedge E_{k,12-l}+\sum_{s=k+1}^{11-k}E_{k,s}\wedge
E_{s,12-k}\right) .
\end{equation*}
The rotation freedom is enclosed in the generators $\widehat{H}_{k}\left(
\gamma \right) =H_{k,n-k+1}+\sum \gamma _{l}H_{l}^{\perp }$. The conditions
(\ref{orth-cart-cond}) define the Cartan elements
\begin{equation*}
H_{k}^{\perp }=\sum_{v=1}^{11}\frac{4k-2}{11}E_{v,v}-\sum_{u=1}^{2k-1}\left(
E_{u,u}+E_{n-u+1,n-u+1}\right) ,
\end{equation*}
This set corresponds to the coordinates $\left\{ \widehat{E}_{k}\right\} $
in which the adjoint action of $\frak{g}_{c}^{\#}$ is reducible,
\begin{equation*}
\begin{array}{l}
\widehat{E}_{1}=E_{2,1}+E_{10,11}, \\
\widehat{E}_{2}=E_{4,3}+E_{8,9}, \\
\widehat{E}_{3}=E_{6,5},
\end{array}
\quad
\begin{array}{l}
\widehat{E}_{4}=E_{8,7}+E_{4,5}, \\
\widehat{E}_{5}=E_{10,9}+E_{2,3}.
\end{array}
\end{equation*}
The pairs $\left\{ H_{k}^{\perp },\widehat{E}_{k}\right\} $ fix the
Jordanian $r$-matrix
\begin{equation*}
r_{J}=\sum_{k=1}^{5}H_{k}^{\perp }\wedge \widehat{E}_{k}.
\end{equation*}
Imposing the CYBE on the sum $r_{rch}\left( \gamma \right) +r_{J}$ ,
\begin{equation*}
\left[ \left[ r_{rch}\left( \gamma \right) +r_{J}\right] \right] =0.
\end{equation*}
we get the necessary and sufficient conditions in the form

\begin{equation*}
\sum_{k=1}^{5}\left[ \widehat{H}_{k}\left( \gamma \right) ,E_{2l,2l-1}\right]
=0;\qquad l=1,\ldots ,5;
\end{equation*}
\begin{equation*}
\begin{array}{l}
\left[ \widehat{H}_{1}\left( \gamma \right) ,E_{10,11}\right] =E_{10,11}, \\
\left[ \widehat{H}_{2}\left( \gamma \right) ,E_{2,3}\right] =E_{2,3}, \\
\left[ \widehat{H}_{3}\left( \gamma \right) ,E_{8,9}\right] =E_{8,9}, \\
\left[ \widehat{H}_{4}\left( \gamma \right) ,E_{4,5}\right] =E_{4,5},
\end{array}
\
\begin{array}{l}
\left[ \widehat{H}_{2}\left( \gamma \right) ,E_{10,11}\right] =-E_{10,11},
\\
\left[ \widehat{H}_{3}\left( \gamma \right) ,E_{2,3}\right] =-E_{2,3}, \\
\left[ \widehat{H}_{4}\left( \gamma \right) ,E_{8,9}\right] =-E_{8,9}, \\
\left[ \widehat{H}_{5}\left( \gamma \right) ,E_{4,5}\right] =-E_{4,5};
\end{array}
\end{equation*}
\begin{eqnarray*}
\left( \widehat{H}_{1}\left( \gamma \right) -\widehat{H}_{2}\left( \gamma
\right) -2H_{1,2}\right) &\sim &H_{1}^{\perp },\left( \widehat{H}_{2}\left(
\gamma \right) -\widehat{H}_{3}\left( \gamma \right) -2H_{9,10}\right) \sim
H_{5}^{\perp }, \\
\left( \widehat{H}_{3}\left( \gamma \right) -\widehat{H}_{4}\left( \gamma
\right) -2H_{3,4}\right) &\sim &H_{2}^{\perp },\left( \widehat{H}_{4}\left(
\gamma \right) -\widehat{H}_{5}\left( \gamma \right) -2H_{7,8}\right) \sim
H_{4}^{\perp }, \\
\left( \widehat{H}_{5}\left( \gamma \right) -2H_{5,6}\right) &\sim
&H_{3}^{\perp };
\end{eqnarray*}
These relations define the unique solution:
\begin{equation*}
\begin{array}{l}
\widehat{H}_{1}=2\sum_{i=1}^{5}\left( H_{2i-1,2i}+H_{i}^{\perp }\right) , \\
\widehat{H}_{2}=2\sum_{i=2}^{5}\left( H_{2i-1,2i}+H_{i}^{\perp }\right) , \\
\widehat{H}_{3}=2\sum_{i=2}^{4}\left( H_{2i-1,2i}+H_{i}^{\perp }\right) ,
\end{array}
\quad
\begin{array}{l}
\widehat{H}_{4}=2\sum_{i=3}^{4}\left( H_{2i-1,2i}+H_{i}^{\perp }\right) , \\
\widehat{H}_{5}=2\left( H_{5,6}+H_{3}^{\perp }\right) .
\end{array}
\end{equation*}
Thus the enlarged chain for $sl(11)$ is found to be
\begin{equation*}
r_{ech}=\sum_{k=1}^{5}\left( \widehat{H}_{k}\wedge
E_{k,12-l}+\sum_{s=k+1}^{11-k}E_{k,s}\wedge E_{s,12-k}\right)
+\sum_{k=1}^{5}H_{k}^{\perp }\wedge \widehat{E}_{k}.
\end{equation*}
This expression is a representative element of the 10-dimensional variety of
solutions

\begin{eqnarray*}
r_{ech} &=&\sum_{k=1}^{5}\xi _{k}\left( \widehat{H}_{k}\wedge
E_{k,12-l}+\sum_{s=k+1}^{11-k}E_{k,s}\wedge E_{s,12-k}\right)
+\sum_{k=1}^{5}\zeta _{k}H_{k}^{\perp }\wedge \widehat{E}_{k}\left( \xi
\right) , \\
\xi _{2,\ldots ,5} &\in &\mathsf{C}^{1}\setminus 0,\xi _{1},\zeta _{k}\in
\mathsf{C}^{1}.
\end{eqnarray*}
The coordinates $\widehat{E}_{k}\left( \xi \right) $ are parameterized as
follows:
\begin{equation*}
\begin{array}{l}
\widehat{E}_{1}\left( \xi \right) =E_{2,1}+\frac{\xi _{1}}{\xi _{2}}%
E_{10,11}, \\
\widehat{E}_{2}\left( \xi \right) =E_{4,3}+\frac{\xi _{3}}{\xi _{4}}E_{8,9},
\\
\widehat{E}_{3}\left( \xi \right) =E_{6,5}.
\end{array}
\quad
\begin{array}{l}
\widehat{E}_{4}\left( \xi \right) =E_{8,7}+\frac{\xi _{4}}{\xi _{5}}E_{4,5},
\\
\widehat{E}_{5}\left( \xi \right) =E_{10,9}+\frac{\xi _{2}}{\xi _{3}}E_{2,3},
\end{array}
\end{equation*}

*************************************************

\end{document}